\documentclass[12pt,oneside]{amsproc}
\usepackage[top=1.25in, bottom=1.25in, left=1.75in, right=1.75in]{geometry}
\usepackage{amsmath,amssymb,amsthm}
\usepackage{enumitem}
\usepackage[dvipsnames]{xcolor}
\usepackage[unicode,breaklinks=true,colorlinks=true]{hyperref}

\numberwithin{equation}{section}
\newtheorem{theorem}{Theorem}[section]

\newtheorem{lemma}[theorem]{Lemma}
\newtheorem{definition}[theorem]{Definition} 

\newtheorem{assumption}[theorem]{Assumption}

\theoremstyle{remark}

\newcommand{\bke}[1]{\left( #1 \right)}
\newcommand{\bkt}[1]{\left[ #1 \right]}

\newcommand{\norm}[1]{\| #1 \|}

\newcommand{\phase}{\phi}
\newcommand{\al}{\alpha}

\newcommand{\de}{\delta}
\newcommand{\e}{\epsilon}

\newcommand{\la}{\lambda}

\newcommand{\si}{\sigma}

\newcommand{\De}{\Delta}
\renewcommand{\th}{\theta}
\newcommand{\Bp}{\dot B_{p,\infty}^{3/p-1}}
\newcommand{\R}{{\mathbb R }}
\newcommand{\RR}{{\mathbb R }}
\newcommand{\N}{{\mathbb N}}
\newcommand{\Z}{{\mathbb Z}}

\newcommand{\pd}{{\partial}}
\newcommand{\nb}{{\nabla}}

\newcommand{\I}{\infty}
\newcommand{\mat}[1]{\begin{bmatrix} #1 \end{bmatrix}}
 \renewcommand{\div}{\mathop{\mathrm{div}}\nolimits}

\newcommand{\EQ}[1]{\begin{equation}\begin{split} #1 \end{split}\end{equation}}
\newcommand{\EQN}[1]{\begin{equation*}\begin{split} #1 \end{split}\end{equation*}}
\DeclareMathOperator*{\esssup}{ess\,sup}
\newcommand{\loc}{\mathrm{loc}} 
\newcommand{\uloc}{\mathrm{uloc}}

\begin{document}

\title[Self-similar solutions to the Navier-Stokes equations]%
{Self-similar solutions to the Navier-Stokes equations: a survey of recent results}
  
\author{Zachary Bradshaw and Tai-Peng Tsai}
\date{\today}
\maketitle 

\begin{abstract}We survey the various constructions of forward self-similar solutions (and generalizations of self-similar solutions) to the Navier-Stokes equations.  We also include and prove an extension of a recent result from \cite{BT3}.
\end{abstract}

\section{Introduction}

The Navier-Stokes equations are a system of partial differential equations that describe the evolution of a viscous incompressible fluid's velocity field $v$ and associated pressure $\pi$.  In three dimensional space they are
\begin{equation} 
\begin{array}{ll}\label{eq:NSE}
 \partial_t v -\Delta v +v\cdot\nabla v+\nabla \pi  = 0
\\[5pt]   \nabla\cdot v = 0
\end{array}\quad \mbox{~in~}\R^3\times [0,\infty),
\end{equation}
and are supplemented with some initial data $v_0$. If the nonlinearity $v\cdot\nabla v$ is omitted, this becomes Stokes system.

Leray proved in \cite{leray} that, if $v_0\in L^2$, then a global in time weak solution $v$ to \eqref{eq:NSE} exists.  Hopf later generalized this to bounded domains (where the problem is supplemented with an appropriate boundary condition) in \cite{hopf}. Leray's construction is based on the a priori bound
\begin{equation}\label{ineq.energy}
\esssup_{0<t'<t}\int |v(x,t')|^2dx + \int_0^t \int 2 |\nabla v(x,t')|^2 dx\,dt'\le \int |v_0(x)|^2dx.
\end{equation}
Formally this is a result of testing \eqref{eq:NSE} against a solution $v$ and noting that the nonlinear term vanishes due to incompressibility.
This energy inequality is identical to that satisfied by solutions to the Stokes system.  Neither Leray nor Hopf were able to say much more  about these weak solutions; this remains true of researchers today.  In particular, we still do not know if Leray's weak solutions are unique or if they are smooth, even with smooth, compactly supported  $v_0$ (partial and conditional results are available, but the general questions remain open). 

Leray noticed that solutions to \eqref{eq:NSE} satisfy a special scaling property: given a solution $v$ of \eqref{eq:NSE}, and $\lambda>0$, it follows that 
\begin{equation}
	v^{\lambda}(x,t)=\lambda v(\lambda x,\lambda^2t),
\end{equation}
is also a solution with associated pressure 
\begin{equation}
	\pi^{\lambda}(x,t)=\lambda^2 \pi(\lambda x,\lambda^2t),
\end{equation}
and initial data 
\begin{equation}
v_0^{\lambda}(x)=\lambda v_0(\lambda x).
\end{equation}
When weak solutions to a system are not clearly regular, irregular behavior can sometimes be found by considering special solutions. Along these lines Leray realized that solutions which are invariant to the above scaling, if defined for negative times, would necessarily develop a singularity at time $t=0$.

We say a solution is self-similar (SS) if it is scaling invariant with respect to the above scaling, i.e.~if $v^\lambda(x,t)=v(x,t)$ for all $\lambda>0$. Hence, Leray's proposed singular solutions would be self-similar.
If this scale invariance holds for a particular $\lambda>1$, then we say $v$ is discretely self-similar with factor $\lambda$ (i.e.~$v$ is $\lambda$-DSS). Similarly $v_0$ can be SS or $\lambda$-DSS.  The class of DSS solutions contains the SS solutions since any SS $v$ is $\lambda$-DSS for any $\lambda>1$. 

Self-similar solutions can be forward (i.e.~defined for $t>0$) or backward (i.e.~defined for $t<0$, like Leray's proposed solutions), and both classes are interesting as sources of irregular behavior such as singularity formation and non-uniqueness.  

For the backward case,
as mentioned above, Leray \cite{leray} proposed the SS solutions as candidates for singularity formation. His original problem on the existence of a nontrivial solution $v(\cdot,t)$ in $W^{1,2}(\R^3) \subset L^3(\R^3)$ for every $t<0$ was excluded in Ne\v cas, {R\accent23 u\v {z}i\v {c}ka}, and {\v Sver\'ak} in \cite{NRS}.   Escauriaza, Seregin, and {\v Sver\'ak} \cite{ESS} gave another proof of the result in \cite{NRS} as a consequence of their $L^\I L^3$ regularity criteria. Tsai proved a localized non-existence result in \cite{Tsai-ARMA} for solutions $v$ satisfying $v(t) \in L^q(\R^3)$, $3<q\le \infty$ for any $t<0$, or $v \in L^{10/3}(B_1 \times (-1,0))$.  These results were recently generalized to Lorentz spaces by Chae and Wolf in \cite{Chae-Wolf2} and by Guevara and Phuc in \cite{GuePhuc}.  Chae and Wolf also proved the non-existence of non-trivial backward \emph{discretely} self-similar solutions \cite{Chae-Wolf3} where the scaling factor $\la$ is close to $1$, (with $N=\sup_{x,t}(x^2-t)^{1/2}|v(x,t)|<\infty$ and 
$\la-1\le \de$ for some sufficiently small $\de>0$ depending on $N$), by reducing the problem to the self-similar case via a contradiction argument.  
There are still several important open problems concerning the existence of backward solutions in larger classes.  For example, solutions that are self-similar modulo a rotation (these will be introduced formally later in this article) do not exist if $u \in L^\infty(-1,0; L^3(\R^3))$ (this follows from the result of \cite{ESS}), but this is not known under the weaker assumption 
\[
|u(x,t)| \le \frac C{|x|+\sqrt{-t}} \quad \text{in }\R^3 \times (-1,0),
\]
even though this assumption does exclude backward self-similar solutions \cite{Tsai-ARMA}.  This is surprising because, at face value, rotated self-similar solutions appear very similar to self-similar solutions (e.g.~they both have stationary ansatzes).  Additionally, backward DSS solutions haven't been ruled out under any condition when $\la$ is significantly larger than one.

For the forward case,
it has been conjectured that {forward} self-similar solutions are a good place to look for non-uniqueness of Leray weak solutions to the Navier-Stokes equations \cite{JiaSverak,JiaSverak2}. 
Self-similar solutions are decided by their behavior at a single time and can therefore be determined by an ansatz which satisfies a stationary system resembling the stationary Navier-Stokes equations.   It is known for certain large data and appropriate forcing, that solutions to the stationary Navier-Stokes boundary value problem are non-unique \cite{Velte,Yudovich,Galdi,Temam}.  In \cite{JiaSverak}, Jia and \v Sver\'ak speculate that similar non-uniqueness results might hold for the stationary profiles of \emph{forward} self-similar solutions, provided the initial data is large in a borderline space such as $L^3_w$. Strong numerical evidence has recently been given by Guillod and \v Sver\'ak  in \cite{GuiSve} supporting the non-uniqueness conjecture in these spaces.
 
Until recently, existence results for forward self-similar and DSS solutions were only known for small data (for small data existence of forward self-similar solutions see \cite{GiMi,Kato,Barraza,CP,Koch-Tataru}).  Such solutions are unique and, therefore, the scale invariance of the solution follows from that of the initial data.
For large data there is no uniqueness result.
 Thus, an existence theory of scale invariant solutions for large data is needed. This article will survey recent advances toward a robust existence theory for forward self-similar and discretely self-similar solutions with large data.  As we will see, the problem is largely solved, but several important outstanding cases remain open. 

This article is structured as follows.  In Section \ref{sec.properties} we will introduce and discuss several classes of solutions to the Navier-Stokes equations with special symmetries (self-similar solutions being an example of such a class).  In Section \ref{sec.classes} we introduce classes of solutions to the Navier-Stokes equations, including ones for initial data that is uniformly locally square integrable or belongs to scaling invariant function spaces.
In Section \ref{sec.overview} we will give an overview of the various constructions of solutions with special symmetries that are presently available.  In Section \ref{sec.construction} we will focus on the construction from \cite{BT1,BT2,BT3} and use it to prove Theorem \ref{thrm.new}, which is an improvement of the existence result in \cite{BT3}.  In particular, Theorem \ref{thrm.new} establishes a local energy inequality for the SS/DSS solution constructed by the method from \cite{BT1,BT3} for initial data in Besov spaces of negative order.

\section{Properties of scaling invariant solutions}\label{sec.properties}

Recall that $v$ is self-similar (SS) if $v(x,t)=\la v(\la x,\la^2 t)$ for all $\lambda>0$ and $v_0$ is SS if $v_0(x)=\la v_0(\la x)$ for all $\la>0$.
If this holds for a particular $\lambda>1$, then $v$ (or $v_0$) is $\lambda$-DSS.  Clearly, if $v$ is SS then it is $\la$-DSS for all $\la>1$. Of course, a solution can also be strictly discretely self-similar.  

Self-similar solutions have a stationary quality in that there exists an ansatz for $v$ in terms of a time-independent profile $u$, namely, 
\begin{equation}\label{ansatz1}
v(x,t) = \frac 1 {\sqrt {t}}\,u\bigg(\frac x {\sqrt{t}}\bigg).
\end{equation} 
The above applies to forward solutions, and we replace $t$ by $-t$ for backward solutions.
The profile $u$ solves the \emph{Leray equations}
\begin{equation} 
\begin{array}{ll}\label{eq:stationaryLeray}
 -\Delta u-\frac 1 2 u-\frac 1 2 y\cdot\nabla u +u\cdot \nabla u +\nabla p = 0&%
\\[5pt]  \nabla\cdot u=0&%
\end{array}
\mbox{~in~}\R^3,
\end{equation}
in the variable $y=x/\sqrt{ t}$.    Note that these equations closely resemble the stationary Navier-Stokes equations.

Discretely self-similar solutions on the other hand are determined by their behavior on time intervals of the form $\la^k\leq t\leq \lambda^{2+k}$.  Thus they have a periodic quality. To make this more precise consider the ansatz
\begin{equation}\label{ansatz2}
v(x,t)=\frac 1 {\sqrt{t}}\, u(y,s),\quad  \pi(x,t) = \frac 1{ t} p(y,s),
\end{equation}
where
\begin{equation}\label{variables2}
y=\frac x {\sqrt{t}},\quad s=\log t.
\end{equation}
Then, the vector field $u$ is $T$-periodic in $s$ with period $T=2\log \lambda$ and solves the \emph{time-dependent Leray equations}
\begin{equation} 
\begin{array}{ll}
\label{eq:timeDependentLeray}
 \partial_s u-\Delta u-\frac 1 2 u-\frac 1 2 y\cdot\nabla u +u\cdot \nabla u +\nabla p = 0& 
\\[5pt]   \nabla\cdot u = 0&
\end{array}
\mbox{~in~}\R^3\times \R,
\end{equation}
where $\nb = \nb_y$ and $\Delta = \Delta_y$. 
Note that the \emph{similarity transform} \eqref{ansatz2}--\eqref{variables2} gives a one-to-one correspondence between solutions to \eqref{eq:NSE} and \eqref{eq:timeDependentLeray}.

A rotational correction can be added to self-similar and discretely self-similar symmetries. For ease of notation,
we will only consider rotations around the $x_3$-axis with matrices
\[
R_s=R( s) = \mat{ 
\cos s & -\sin s & 0 \\
\sin s & \cos s & 0 \\
0 & 0 & 1 }.\]
Note $R( s)R(\tau)=R(\tau)R( s)$ for any $s,\tau \in \R$, and
\[
\frac d{ds} R(s) = JR(s) = R(s)J, \quad 
J= \mat{ 
0 & -1 & 0 \\
1&0 & 0 \\
0 & 0 & 0 }.\]

A vector field $v(x,t)$ is said to be \emph{rotated self-similar} (RSS) if, for some fixed $\alpha\in \R$ and for \emph{all} $\la>0$, $x\in \R^3$, and $t>0$,
\begin{equation}
\label{v-RSS}
v(x,t) = \la R(- 2\alpha \log \la)\, v\! \bke{\la R(2 \alpha \log \la) x, \la^2 t}.
\end{equation}
The constant $\al$ will be called the \emph{angular speed}, and is understood relative to the new time variable $s$ to be defined in \eqref{eq.v.to.u}.
An RSS vector field is always DSS with any factor $\la>1 $ such that $2 \alpha \log \la \in 2\pi \mathbb{Z}$.
When $\alpha=0$ it becomes SS.  Hence SS $\subsetneq$ RSS $\subsetneq$ DSS. 

Setting $\la=t^{-1/2}$, any RSS vector field $v$ satisfies
\EQ{
v(x,t) = R(\al \log {t}) \frac 1{\sqrt{t}}\, v\!\bke{R(-\al \log {t}) \frac x{\sqrt{t}} , 1},
} 
for all $x\in\R^3$ and $t>0$.
Thus the value of $v$ is determined by its value at any fixed time, and given any profile at a fixed time we can construct an RSS vector field.  Hence RSS solutions have the same stationary quality as SS solutions.

A vector field $v(x,t)$ is said to be \emph{rotated discretely self-similar} (RDSS) if, for some $\la>1$ (not necessarily all $\la>1$) 
and some $\phase \in \R$, 
\begin{equation}
\label{v-RDSS}
v(x,t) = \la R(-\phase)\, v\! \bke{\la R(\phase) x, \la^2 t},
\end{equation}
for all $x\in\R^3$ and $t>0$.
We call $\la$ the \emph{factor} and $\phase$ the \emph{phase}.
When $\phase\in 2\pi \mathbb{Z}$ we recover $\la$-DSS vector fields. 
If $n\phase=2\pi m$ for some integers $n>0$ and $m$, then $v$ is DSS with factor $\la^n$. If $\frac\phase{2\pi}$ is irrational, in general $v$ is not DSS. 
For any $t>0$ let $\tau(t)\in [1,\lambda^2)$ satisfy $\tau=\lambda^{2k}t$ for some $k\in \Z$.  Then,
\[
v(x,t)= \la^k R(-k\phase)\, v\! \bke{\la^k R(k \phase) x, \tau},
\]
i.e., $v$ is decided entirely by its values on $t\in [1,\lambda^2)$.
{Note that an RSS vector field with angular speed $\al$ is always RDSS for any factor $\la>1$ with phase $\phase=2\al \log \la $.}

In summary, the inclusions between these classes are
\[
\mbox{SS}\subsetneq \mbox{RSS} \subsetneq \mbox{DSS} \subsetneq \mbox{RDSS}.
\] 

For initial data, we say a vector field
$v_0(x):\R^3\to \R^3$ is RSS if, for some $\alpha\in \R$ and all $x\in\R^3$ and $\la>0$,
\begin{equation}
\label{v0-RSS}
v_0(x) = \la R(- 2 \alpha \log \la) v_0\bke{\la R( 2\alpha \log \la) x}, 
\end{equation}
and is
RDSS if for some $\la>1$ and some $\phase \in \R$, 
\begin{equation}
\label{v0-RDSS}
v_0(x) = \la R(-\phase) v_0\bke{\la R(\phase) x}, 
\end{equation}
for all $x\in \R^3$.
Note that, like SS initial data, if $v_0$ is RSS then it is determined by its values on the unit sphere. Similarly, if $v_0$ is RDSS, then it is determined by its values on $\{ x: 1\leq |x|<\la \}$.

RSS and RDSS solutions also have ansatzes which are stationary or time-periodic.  Indeed, let
\begin{equation}
\label{eq.v.to.u} 
v(x,t) =\frac 1{\sqrt t} R_\theta u(y,s), \quad  \pi(x,t) = \frac 1{ t} p(y,s),\quad y = R_\theta^T \frac x {\sqrt t},\quad s=\log t,
\end{equation}
for some function $\th(s)$. Then $u$ satisfies  
\EQ{\label{eq:Leray-rot}
  \partial_s u  + \dot \th J u -  (\dot \th J y)\cdot \nb u
  &-  \frac 12 u -  \frac y2\cdot\nabla u -\Delta u +u\cdot\nabla u+\nabla p  = 0,
 \\ & \nabla \cdot u = 0,
}
in $\R^3\times \R$.
If $v(x,t)$ is an RSS solution of \eqref{eq:NSE} satisfying \eqref{v-RSS}, then $u(y,s)$ is a stationary solution of \eqref{eq:Leray-rot} with constant $\dot \theta = \alpha$, if we take  $\theta(s)=\al s$.
For any RDSS solution $v(x,t)$ of \eqref{eq:NSE} satisfying \eqref{v-RDSS} with factor $\lambda>1$ and phase $\phase$, 
let 
\begin{equation}
\label{alpha-choice}
T=2 \log \la,\quad
 \al_k= \frac {2k\pi + \phase}T ,
\end{equation}
for an arbitrary integer $k\in \mathbb{Z}$. Then $v(x,t)$ 
corresponds to a periodic solution $u(y,s)$ of \eqref{eq:Leray-rot} with constant $\dot \theta = \alpha_k$ and period $T$. To be definite we will take $\al=\al_0 = \frac \phase T$.

\section{Solution classes}\label{sec.classes}

In this section we introduce various function spaces and classes of solutions to the Navier-Stokes equations.  These will be highly relevant to our subsequent discussion of self-similar solutions and their generalizations.

\subsection{Weak solutions}

As mentioned in the introduction, Leray constructed weak solutions to the Navier-Stokes equations in \cite{leray}.  His solutions satisfied the global energy inequality \eqref{ineq.energy}.  In honor of Leray's contribution, any weak solution with data in $L^2$ satisfying \eqref{ineq.energy} is called a \emph{Leray weak solution}.  

Lemari\'e-Rieusset introduced a local analogue of Leray's solutions in \cite{LR} called \emph{local Leray solutions}.  We recall the definition in full.  For $q \in [1,\infty)$, we say $f \in L^q_\uloc$ if $f$ is defined in $\R^3$ and
\[
\norm{f}_{L^q_\uloc} = \sup_{x_0\in \R^3} \norm{f}_{L^q(B_1(x_0))}<\infty.
\]

\begin{definition}[Local Leray solutions]\label{def:localLeray} A vector field $v\in L^2_{\loc}(\R^3\times [0,\infty))$ is a local Leray solution to \eqref{eq:NSE} with divergence free initial data $v_0\in L^2_{\uloc}$ (and zero force) if:
\begin{enumerate}[label=$($\alph*$)$]
\item for some $\pi\in L^{3/2}_{\loc}(\R^3\times [0,\infty))$, the pair $(v,\pi)$ is a distributional solution to \eqref{eq:NSE},
\item for any $R>0$, $v$ satisfies
\begin{equation}\notag
\esssup_{0\leq t<R^2}\,\sup_{x_0\in \R^3}\, \int_{B_R(x_0 )} |v(x,t)|^2\,dx + \sup_{x_0\in \R^3}\int_0^{R^2}\! \!\int_{B_R(x_0)} |\nabla v(x,t)|^2\,dx \,dt<\infty,\end{equation}
\item for any $R>0$, $v$ satisfies
\begin{equation}\notag 
\lim_{|x_0|\to \infty} \int_0^{R^2}\!\! \int_{B_R(x_0 )} | v(x,t)|^2\,dx \,dt=0,
\end{equation}
\item for all compact subsets $K$ of $\R^3$ we have $v(t)\to v_0$ in $L^2(K)$ as $t\to 0^+$,
\item $v$ is suitable in the sense of Caffarelli-Kohn-Nirenberg, i.e., for all cylinders $Q$ compactly supported in  $ \R^3\times(0,\infty )$ and all non-negative $\phi\in C_0^\infty (Q)$, we have 
\EQ{\label{ineq.localv}
&%
2\iint |\nabla v|^2\phi\,dx\,dt 
\\&\leq 
\iint |v|^2(\partial_t \phi + \Delta\phi )\,dx\,dt +\iint (|v|^2+2\pi)(v\cdot \nabla\phi)\,dx\,dt.
}
\end{enumerate}
\end{definition}

In \cite{LR} Lemari\'e-Rieusset constructed global in time  local Leray solutions if $v_0$ belongs to $E_2$, the closure of $C_0^\infty$ in the $L^2_{\uloc}(\R^3)$ norm.  See Kikuchi-Seregin \cite{KiSe} for more details. In particular, condition 3 justifies a formula of the pressure $\pi$ in terms of the velocity $v$, see \cite[(1.9)]{KiSe} and \cite[(3.3)]{JiaSverak}.

Local Leray solutions are known to satisfy a useful a priori bound.  Let $\mathcal N (v_0)$ denote the set of all local Leray solutions with initial data $v_0$. The following estimate is known for local Leray solutions (see \cite{JiaSverak}): for all $\tilde v\in \mathcal N (v_0)$ and $r>0$ we have
\begin{equation}\label{ineq.apriorilocal}
\esssup_{0\leq t \leq \sigma r^2}\sup_{x_0\in \RR^3} \int_{B_r(x_0)}\frac {|\tilde v|^2} 2 \,dx\,dt + \sup_{x_0\in \RR^3}	\int_0^{\sigma r^2}\!\! \int_{B_r(x_0)} |\nabla \tilde v|^2\,dx\,dt <C A _r,
\end{equation}
where 
\begin{equation} 
\label{ineq.apriorilocalconstant}
A_r= \sup_{x_0\in \RR^3} \int_{B_r(x_0)} \frac {|v_0|^2}2\,dx, \quad
\sigma(r) =c_0\, \min\big\{r^2A_r^{-2} , 1  \big\},
\end{equation}
for a small universal constant $c_0$.  

Local Leray solutions are currently only defined on the whole space because the pressure formula breaks down on domains possessing boundaries. 

\subsection{Mild solutions}

There is a rich literature about the global well-posedness of solutions to the Navier-Stokes equations when the data belongs to a scaling invariant function space (see \cite{LR} for a review of this).  Note that $X$ is scaling invariant if $\|u\|_X= \|u^{\la} \|_X$ for any $\la$.  Since we are working over $\R^3$, the most obvious example of a scaling invariant function space is $L^3$.  Global well-posedness is known in $L^3$ for small data. But, $L^3$ does not include non-trivial self-similar data because any such data will look like $|x|^{-1}$ on rays emanating from the origin.

The weak Lebesgue space $L^3_w$ is slightly larger than $L^3$ and includes $|x|^{-1}$, making it a natural space to study self-similar initial data.  Recall that $f\in L^3_w(\R^3)$ if and only if
$\|f\|_{L^3_w(\R^3)}<\I$,
where 
\begin{equation}\label{weakL3}
\|f\|_{L^3_w(\R^3)} = \sup_{s>0}\,s \,m(f,s)^{1/3},
\end{equation}
and $m(f,s)$ is the distribution function of $f$ given by
\[
m(f,s)=|\{ x\in \R^3: | 	f(x)|>s  \}|.
\]  
Here $|S|$ denotes the three dimensional Lebesgue measure of a given set $S$.
Let $L^3_{w,\si}(\R^3)$ be the subspace of $L^3_w$ of divergence free vector fields.   
In $L^3_w$ the best global well-posedness result is the following: if $u_0\in L^3_{w,\si}$ and $\|u_0\|_{L^3_w}<\epsilon_0$ for a universal, small constant $\e_0$, then there exists a global in time \emph{mild} solution $u$ which is also a strong (i.e.~classical) solution and is unique among solutions belonging to $L^\I(0,\I;L^{3}_w)$ \cite{Barraza,LR}.  By mild solution we mean a solution to  \eqref{eq:NSE} satisfying an integral formulation.  Such solutions are nonlinear perturbations of $e^{t\Delta}v_0$.

Note that $L^3_w(\RR^3)$ embeds continuously into the space of uniformly locally square integrable functions $L^2_{\uloc}(\RR^3)$.  Thus, for any initial data in $L^3_{w,\si}$, one may construct a global-in-time {local Leray solutions} in the whole space. However, this solution 
may not be unique for large data, and hence is not necessarily self-similar.

Analogous wellposedness results for small data can be formulated in the Besov spaces $\Bp$ ($1<p<\I$) and the Koch-Tataru space $BMO^{-1}$, which we now introduce.  Besov spaces can be defined using the Littlewood-Paley decomposition.  Fix an inverse length scale $\la>1$.  
Let $B_r$ denote the ball of radius $r$ centered at the origin in $\R^3$.  Fix a non-negative, radial cut-off function $\chi\in C_0^\infty(B_{1})$ so that $\chi(\xi)=1$ for all $\xi\in B_{1/\la}$. Let $\phi(\xi)=\chi(\lambda^{-1}\xi)-\chi(\xi)$ and $\phi_j(\xi)=\phi(\lambda^{-j}\xi)$.  For a vector field $u$ of tempered distribution, let $\Delta_j u=(\mathcal F^{-1}\phi_j)*u$ for $j\in \N_0 $ and $\Delta_{-1}=(\mathcal F^{-1}\chi)*u$. Then, $u$ can be written as\[u=\sum_{j\geq -1}\Delta_j u.\]
If $(\mathcal F^{-1}\chi(\la^{-j}\cdot))*u\to 0$ as $j\to -\infty$ in the space of tempered distributions, then for $j\in \Z$ we define $\dot \Delta_j u = \mathcal F^{-1}\phi_j*u$ and have
\[u=\sum_{j\in \Z}\dot \Delta_j u.\] 
For $s\in \R$, $1\leq p,q\leq \infty$, the non-homogeneous Besov spaces include tempered distributions modulo polynomials for which the norm
\begin{align*}
&\|u\|_{B^s_{p,q}}:= 
\begin{cases} 
 \bigg(\sum_{ j\geq -1} \big(    \lambda^{sj} \|\Delta_j u \|_{L^p(\R^n)}  \big)^q \bigg)^{1/q}   & \text{ if } q<\infty  
\\ \sup_{j\geq -1} \lambda^{sj} \| \Delta_j u \|_{L^p(\R^n)} & \text{ if } q=\infty
\end{cases}, 
\end{align*}is finite, while the homogeneous Besov spaces include tempered distributions modulo polynomials for which the norm 
\begin{align*}
&\|u\|_{\dot B^s_{p,q}}:= 
\begin{cases} 
 \bigg(\sum_{ j\in \Z} \big(    \lambda^{sj} \|\dot \Delta_j u \|_{L^p(\R^n)}  \big)^q \bigg)^{1/q}   & \text{ if } q<\infty  
\\ \sup_{j\in \Z} \lambda^{sj} \| \dot \Delta_j u \|_{L^p(\R^n)} & \text{ if } q=\infty
\end{cases},
\end{align*}
is finite.
Note that the partition of unity can be taken to be $\la$-adic for any $\la>1$.  The resulting Besov norm is equivalent to any dyadic norm \cite{BT3}.
For $3\leq p<\I$, a small data global well-posedness theory exists for data in $\Bp$ that is identical to that for data in $L^3_w$ \cite{BCD} (with $L^3_w$ replaced by $\Bp$).  Since $L^3_w\subsetneq \Bp$ whenever $p>3$, this generalizes the result of \cite{Barraza}.

The Koch-Tataru space $BMO^{-1}$ is the largest critical function space in which the Navier-Stokes equations are globally well-posed for small data and consists of distributions which are derivatives of elements of $BMO$, the class of functions with bounded mean oscillation.  The existence of global in time mild strong solutions was proved by Koch and Tataru in \cite{Koch-Tataru}.
Note that $\Bp\subsetneq BMO^{-1} \subsetneq \dot B_{\I,\I}^{-1}$.  We have observed that small data well-posedness holds for the first two spaces; it fails for the third:  In \cite{BoPa}, Bourgain and Pavlovic proved the last point by showing that, for any time length $\delta$, any $\e>0$, and any $M>0$, there exists an initial data $u_0$ so that $\|u_0\|_{\dot B_{\I,\I}^{-1}}<\e$ and $\|u(\delta)\|_{\dot B_{\I,\I}^{-1}}>M$ (this is called norm inflation).

\section{Existence results} \label{sec.overview}

\subsection{Strong self-similar solutions for small data}

In the previous section we mentioned global-in-time well posedness is known for small data in certain critical spaces including $L^3_w$, $\Bp$, and $BMO^{-1}$.  These results give the existence of SS/RSS/DSS/RDSS solutions as a corollary.  Indeed, assume that $v_0$ is SS and $\|v_0\|_X$ is small where $X$ is any of the spaces just mentioned.  Then, there exists a strong mild solution $v$ to \eqref{eq:NSE} that belongs to and is unique in $C_w([0,\I);X)$.
  But, for every $\la>0$, $v^{\la }$ is also a solution and belongs to $C_w([0,\I);X)$. Since $v_0$ is SS, $v_0^{\la }=v_0$ and, by uniqueness in $C_w([0,\I);X)$, it follows that $v=v^{\la}$ for all $\la$, i.e.~$v$ is self-similar.  The same argument goes through for RSS, DSS, and RDSS data. 
   
For small data, self-similar solutions were studied in the eighties and nineties by \cite{GiMi, Kato, Barraza, CP, Koch-Tataru}.  
Giga-Miyakawa \cite{GiMi} constructed solutions to the vorticity equations assuming the initial vorticity is in a critical space of measures.
Kato gave the construction in \cite{Kato} assuming the initial velocity was small in a critical Morrey space. 
 In 1996, Barraza \cite{Barraza} constructed solutions for small data in $L^3_w$, as did Cannone and Planchon \cite{CP} for small data in $\Bp$ where $p>3$.  Finally the solution is constructed in $BMO^{-1}$ by Koch-Tataru \cite{Koch-Tataru}.  
The regularity properties of Koch-Tataru  solutions are studied by Miura and Sawada \cite{MiuraSawada} and by
 Germain, Pavlovi\'c, and Staffilani \cite{GePaSt}.
Decay and asymptotic properties of small self-similar solutions have been further examined by Brandolese in \cite{Bran}. 

\subsection{Strong self-similar solutions for large, smooth data}

As mentioned above, global well-posedness is only known in $L^3_w$ for small data.  For large data, only weak solutions are known to exist.  Indeed, since $L^3_w$ is in the closure of $C_0^\I$ under $L^2_{\uloc}$-norm, for any $v_0\in L^3_w$, there exists a global in time local Leray solution $v$ \cite{KiSe,LR}.  Barker and Seregin have recently given another approach to constructing global solutions for large data in $L^3_w$ \cite{BaSe}.  Because uniqueness is not known in these classes, it is unclear if the solutions constructed in \cite{KiSe,LR,BaSe} inherit the scaling properties of the initial data.  Thus, for large data in $L^3_w$, we do not get the existence of self-similar solutions for free by these methods.  

It was thus surprising when, in 2014, Jia and \v Sver\'ak  constructed a forward self-similar solution for \emph{large $-1$-homogeneous} initial data which is \emph{locally H\"older continuous} away from the origin \cite{JiaSverak}.  Their proof used Leray-Schauder degree theory, the main ingredient of which are \emph{a priori} bounds for self-similar solutions, existence and uniqueness for small data, and global compactness.  Most of the work lies in establishing H\"older estimates for the solutions, which depend on the fact that the solutions live in the local Leray class.  

In \cite{Tsai-DSSI}, Tsai gave a similar result for $\lambda$-DSS solutions with factor \emph{close to one} where closeness is determined by the local H\"older norm of $v_0$ away from the origin.  It is also shown in \cite{Tsai-DSSI} that the closeness condition on $\lambda$ can be eliminated if the initial data is axisymmetric with no swirl. The approach is similar to \cite{JiaSverak}.
In \cite{KT-SSHS}, Korobkov and Tsai constructed self-similar solutions on the half space (their approach also works on the whole space) for appropriately smooth initial data.  The approach here differs from \cite{JiaSverak} and \cite{Tsai-DSSI} in that the existence of a solution to the stationary Leray equations \eqref{eq:stationaryLeray} is established directly.  It also gives a second proof of the main result of \cite{JiaSverak}.  {A new approach is necessary in 
\cite{KT-SSHS} due to lack of spatial decay estimates, which gives  \emph{global compactness} needed for the Leray-Schauder theorem in 
\cite{JiaSverak} and \cite{Tsai-DSSI}.}

\subsection{Weak self-similar solutions for large, rough data}

The constructions in \cite{JiaSverak,Tsai-DSSI,KT-SSHS} necessarily produce strong solutions; this is why they require the initial data to be H\"older continuous.  If we instead  use construction methods which only necessarily yield weak solutions, then we can weaken the assumptions imposed on the initial data.  Indeed, this idea is the main motivation for the sequence of papers by the authors, \cite{BT1}, \cite{BT2}, and \cite{BT3}, where SS/RSS/DSS/RDSS solutions are constructed for large, possibly rough data.  In \cite{BT1} and \cite{BT2} we considered data in $L^3_w$ and arrived at the following theorem (which is an amalgam of the results of \cite{BT1,BT2}). 

\begin{theorem}\label{thrm:BT1}
Let $v_0$ be a SS vector field and belong to $L^3_{w,\si}(\R^3)$. Then, there exists a local Leray solution $v$ to \eqref{eq:NSE} which is SS and additionally satisfies
\begin{equation}\notag
 \|  v(t)-e^{t\Delta}v_0 \|_{L^2(\R^3)}\leq C_0\,t^{1/4}
\end{equation}
for any $t\in (0,\infty)$ and some constant $C_0=C_0(v_0)$.

The above statement is also true with ``SS'' replaced by any of the following: 
\begin{itemize}
\item ``RSS for any given $\al\in \R$'',
\item  ``DSS for any given $\la>1$'',
\item ``RDSS for any given $\al\in \R$ and $\la>1$''.
\end{itemize}
\end{theorem}

This result is more general than all earlier results.  In particular, the initial data can be discontinuous or even singular away from the origin.  Furthermore, it is valid for any choice of available parameters.  The price paid for this generality is a loss of \emph{a priori} regularity of the solutions; the solutions from \cite{BT1} are not guaranteed to be smooth by their construction.  In the SS and RSS cases, smoothness follows after the fact using Gruji\'c's result which implies that any SS local Leray solution is smooth \cite{Grujic} (that proof also works in the RSS case).  In the DSS and RDSS cases, smoothness is not known generally, but is expected for values of $\la$ close to $1$.

An analogous result can be formulated on the half space, as is shown in \cite{BT2}.  Adapting the proof of \cite{BT1} to that case requires several technical modifications. 
In particular, the solution class is weakened, and we do not impose the local energy inequality. 
 In comparison, self-similar solutions for smooth data on $\R^3_+$ were constructed in \cite{KT-SSHS} but strictly discretely self-similar solutions had not been constructed prior to \cite{BT2}.

Because existence of SS solutions is known for small data in spaces larger than $L^3_w$, in particular the Besov spaces $\Bp$, it is natural to try to extend Theorem \ref{thrm:BT1} in that direction.   Partial progress was made in \cite{BT3} to the Besov spaces $\Bp$ for $3<p<6$ for SS and DSS data. 
The following theorem is a slight refinement of the main result of \cite{BT3}.  

\begin{theorem}\label{thrm.new} Fix $p\in (3,6)$.  Assume $v_0:\R^3\to \R^3$ is divergence free, belongs to $\Bp$, and is self-similar. Then there exists a self-similar distributional solution $v$ and pressure distribution $\pi$ to 3D NSE on $\R^3\times (0,\infty)$.  
Furthermore, $v$ and $v_0$ can be decomposed as $a+b$ and $a_0+b_0$ respectively so that $a_0\in L^3_w$, $b_0$ is small in $\Bp$, and $a$ and $b$ satisfy the following properties:
\begin{itemize}
\item $a$ and $b$ are self-similar,
\item $b$ is a global strong solution to \eqref{eq:NSE} for initial data $b_0$ having associated pressure $p_b$,
\item there exists a distribution $p_a$ so that $a$ and $p_a$ solve 
\EQ{\label{eq.a}
\partial_t a- \Delta a + a\cdot\nabla a+ a\cdot\nabla b + b\cdot \nabla a +\nabla p_{a}=0,
\quad \div a=0,
}
in the sense of distributions,
\item $a(t)$ converges to $a_0$ in the sense that
\EQ{\label{th1.3-1}
\|a(t)-e^{t\Delta}a_0\|_{L^2} &\leq C_2 t^{1/4},\\
\int_{0}^t \norm{a(\tau)-e^{\tau\Delta}a_0}_{L^r}^q d\tau  &\le C_r t^{q/4},\quad \forall\, r \in (2,6],
}
for some constant $C_r(v_0)$ with $\frac 3r + \frac 2q=\frac 32$,
\item  $b(t)$ converges to $b_0$ in the sense that
\begin{equation}\label{th1.3-2}
\norm{b(t)-e^{t\Delta}b_0}_{L^r} \le C_r \norm{b_0}_{\Bp}^2 t^{-\frac 12 + \frac 3{2r}}, \quad \forall\, r \in 
\left[\frac p2, \frac {3p}{6-p}\right),
\end{equation}
for some constant $C_r$, 
\item  $v=a+b$ and $\pi = p_a + p_b$ satisfy the local energy inequality \eqref{ineq.localv}
for any  
non-negative $\phi\in C^\I_0(\R^3\times \R^3_+)$.

\end{itemize}

The above statement is also true with ``self-similar'' replaced by {``DSS for any given $\la>1$''.} 
\end{theorem}

We will prove Theorem \ref{thrm.new} in Section \ref{sec.construction}.

\medskip

 \noindent\emph{Comments on Theorem \ref{thrm.new}}:
\begin{enumerate}[label=(\alph*)]
\item Theorem \ref{thrm.new} improves  the main result of \cite{BT3} in that 
it includes the local energy inequality \eqref{ineq.localv}.
 It is however restricted to the whole space setting, not for $\R^3_+$.
\item If $b_0$ is smooth, for $\phi \in C^\infty_0(\R^3 \times [0,\infty))$, possibly non-zero at $t=0$, 
we can replace $\phi$ in \eqref{ineq.localv}
by $\phi \theta(t/\e)$, where $\theta(t)=1$ for $t>1$ and $\theta(t)=0$ for $t<1/2$, and take limits $\e \to 0_+$ to get the usual form of local energy inequality %
\EQ{ \label{ineq.localv-2}
& 2\iint |\nabla v|^2\phi\,dx\,dt 
\leq \int |v_0|^2 \phi|_{t=0} dx  
\\
&\quad+ \iint |v|^2(\partial_t \phi + \Delta\phi )\,dx\,dt
 +\iint (|v|^2+2p)(v\cdot \nabla\phi)\,dx\,dt.
}
For general $b_0 \in \Bp$, however, $b_0$ may not be in $L^2_\loc$ (see \cite[\S6]{BT3}), hence \eqref{ineq.localv-2} is not meaningful.
\end{enumerate}

Theorem \ref{thrm:BT1} has also been generalized in a different direction by Lemari\'e-Rieusset in \cite{LR2} and by Chae and Wolf in \cite{Chae-Wolf}, both working on the whole space.  In \cite{LR2}, Lemari\'e-Rieusset uses the Leray-Schauder approach to first construct self-similar solutions for initial data $v_0$ satisfying $|v_0(x)|\lesssim |x|^{-1}$.  This construction is more general than that in \cite{JiaSverak} but less general than that in \cite{BT1}.  But, Lemari\'e-Rieusset also noticed that provided $v_0$ is self-similar, then $v_0\in L^2_{\loc}$ if and only if $v_0\in L^2_{\uloc}$.  And, furthermore, if $v_0$ is self-similar and belongs to $L^2_{\uloc}$, then it can be approximated by a sequence $v_0^{(k)}$ where each $|v_0^{(k)}(x)|\lesssim |x|^{-1}$. Then, his first construction gives local Leray solutions for each $v_0^{(k)}$ and, because local Leray solutions satisfy the a priori bound \eqref{ineq.apriorilocal} depending only on the $L^2_{\uloc}$ norm of their initial data, these will converge to a SS local Leray solution with $L^2_{\loc}$ data.  Thus, Lemari\'e-Rieusset is able to prove in \cite{LR2} that, for any initial data $v_0\in L^2_{\loc}$ that is incompressible and self-similar, there exists a self-similar local Leray solution evolving from $v_0$. This argument breaks down for DSS solutions since $L^2_{\loc}\cap DSS \neq L^2_{\uloc}\cap DSS$.  Chae and Wolf, on the other hand, use an entirely new method to construct $\la$-DSS solutions for any $\la>1$ and initial data $v_0\in L^2_{\loc}(\R^3)$. These solutions live in the class of ``local Leray solutions with projected pressure,'' which means they satisfy a modified local energy inequality instead of the classical local energy inequality of \cite{CKN}.  We expect this result can be refined slightly; in particular, it seems possible to construct distributional solutions to \eqref{eq:NSE} for any SS/DSS data (or RSS/RDSS data) which satisfy the local energy inequality in the sense of \cite{CKN}. This will be addressed in the future.

\subsection{Next directions}

Theorem \ref{thrm.new} is the strongest available in the scale of critical spaces at the time of writing this survey.  We expect that solutions exist for SS/DSS data in $\Bp$ for any $3<p<\I$, that is, the assumption that $p<6$ is likely not needed.  Indeed, using a lemma from \cite{BT3} (see Lemma \ref{lemma.profileslicing} in the next section), it is possible to show that any DSS datum in $\Bp$ can be approximated by data in $L^3_w$. Each approximate data then gives rise to a DSS solution to \eqref{eq:NSE} by Theorem \ref{thrm:BT1}.  It seems that these solutions \emph{should} converge to a DSS solution to \eqref{eq:NSE} provided some presently unknown a priori bound holds for data in spaces larger than $L^2_{\uloc}$, but we have not found a proof of this.

Although we constructed RSS/RDSS solutions for data in $L^3_w$, 
our method for SS/DSS data in $\Bp$ does not extend to RSS/RDSS data in $\Bp$.
In particular, we have not been able to prove a counterpart of Lemma \ref{lemma.profileslicing} for such data,
since rotations do not change the size but do change the frequency.
It would be interesting to decide if such a construction is intrinsically possible.

Recall that for small data existence is also known in $BMO^{-1}$.  The $\Bp$ argument outlined above has no hope of working in $BMO^{-1}$ because, it appears, Lemma \ref{lemma.profileslicing} breaks down for $v_0\in BMO^{-1}$.  This is because $BMO^{-1}$ is an $L^\I$ based space while $\Bp$ is an $L^p$ based space ($p<\I$) -- i.e.~$\|\cdot \|_{BMO^{-1}}$ is computed by taking the supremum over quantities computed on dyadic cubes of (and above) a fixed scale while $\|\cdot \|_{\Bp}$ is computed by taking the $l^p$ norm over quantities computed on dyadic cubes of a fixed scale (at least, in the wavelet characterization of norms).  In $l^p$, the tail of a convergent series can be made small if $p<\I$ -- this is the idea behind Lemma \ref{lemma.profileslicing} and clearly fails in an $l^\I$ based space.
Thus we expect that a totally new approach is needed to address the existence of large data SS/RSS/DSS/RDSS solutions for data in $BMO^{-1}$.

\section{Constructing DSS/SS solutions: a case study}\label{sec.construction}

In this section we will follow the procedure of \cite{BT3} to construct forward SS/DSS solutions with data in the Besov spaces $\Bp(\R^3)$ where $3<p<6$.   In comparison to \cite{BT3} we include an additional step that reveals more about the structure of the solutions -- in particular, the solution satisfies a local energy inequality away from $t=0$.

We will use the following function spaces:
\begin{align*}
&\mathcal V=\{f\in C_0^\infty({ \R^3;\R^3}) ,\, \nabla \cdot f=0 \},
\\& X = \mbox{the closure of~$\mathcal V$~in~$H^1(\R^3)$} ,
\\& H = \mbox{the closure of~$\mathcal V$~in~$L^2(\R^3)$}.
\end{align*}
Let $X^*(\R^3)$ denote the dual space of $X(\R^3)$. 
Let $(\cdot,\cdot)$ be the $L^2(\R^3)$ inner product and $\langle\cdot,\cdot\rangle$ be the dual product for $H^1$ and its dual space $H^{-1}$, or that for $X$ and $X^*$.
Denote by  $\mathcal D_T$ the collection of all smooth divergence free vector fields in $\R^3 \times \R$ which 
are time periodic with period $T$ and whose supports are compact in space.  

We treat the DSS and SS cases separately, starting with the former.  

\subsection{Proof of Theorem \ref{thrm.new}, DSS case}

Elements of this proof are identical to those in \cite{BT3} and some details are omitted here.  The local energy estimate \eqref{ineq.localv} was not considered in \cite{BT3} and is therefore our main focus.  

\medskip
\noindent{\bf Step 1: formulation of an auxiliary problem.}  Recall that if $v$ is a DSS solution for some initial data $v_0$, then $u$ defined by 
\eqref{ansatz2}-\eqref{variables2}
satisfies \eqref{eq:timeDependentLeray}, i.e.
\EQN{
  \partial_s u  
  -  \frac 12 u -  \frac y2\cdot\nabla u - &\Delta_y u +u\cdot\nabla u+\nabla p  = 0,
 \\ & \nabla \cdot u = 0,
}
and, furthermore, $u$ is time periodic with period $T=2\log \la$. 
Therefore, it suffices to construct $u$ and then show the corresponding vector field $v$ satisfies the properties in the statement of Theorem \ref{thrm.new}. 

Let $\epsilon_0>0$ be small enough that there exists a global strong solution to the Navier-Stokes equations for any data in $\Bp$ with norm smaller than $\e_0$.   In \cite{BT3} we proved the following lemma. 

\begin{lemma}[{\cite[Lemmas 2.2, 5.2]{BT3}}]
\label{lemma.profileslicing}
Let $f$ be a $\la$-DSS, divergence free vector field in $\R^3$, and belong to $\dot B_{p,\infty}^{3/p-1}$ for some $\la \in (1,\infty)$ and $p\in (3,\infty)$. For any $\epsilon>0$, there exist divergence free $\la$-DSS distributions  $a\in L^3_w$ and $b\in \Bp$ so that 
 $f = a + b$ and $ \|b  \|_{\dot B_{p,\infty}^{3/p-1}}<\epsilon$.
{If $f$ is self-similar, then so are $a$ and $b$.} 
\end{lemma}

Apply Lemma \ref{lemma.profileslicing} for $v_0$ and the prescribed value $\epsilon_0$.  Then, $v_0=a_0+b_0$ where $b_0\in \Bp$, $\|b_0\|_{\Bp}<\e_0$, and $a_0\in L^3_w$.  Let $b$ be the strong mild solution to the Navier-Stokes equations for the initial data $b_0$ and let $p_b$ denote its pressure.  By uniqueness in the class of strong mild solutions, $b$ is DSS.  Therefore, there exists an ansatz $B$ for $b$ which is time-periodic with period $T=2\log \la$, divergence free, and which satisfies
\EQ{\notag 
  \partial_s B  
  &-  \frac 12 B -  \frac y2\cdot\nabla B -\Delta_y B +B\cdot\nabla B+\nabla p_B  = 0.
}

Let $A=u-B$. Then, $A$  satisfies   
\EQ{ 
\label{eq:wholeSpaceLeray}
  \partial_s A  
  -  \frac 12 A -  \frac y2\cdot\nabla A -\Delta_y A
   +&A\cdot\nabla A + A\cdot\nabla B + B\cdot \nabla A+\nabla p_A  = 0,
   \\ & \div A=0.
}
Letting $a(x,t)$ be the image of $A(y,s)$ under 
\eqref{ansatz2}-\eqref{variables2}, we see that $a$ satisfies \eqref{eq.a}
and (formally) converges to $a_0$ as $t\to 0^+$.  Thus, constructing $A$ (in some specific context) will, after changing variables, give us the desired solution $v=a+b$ to the Navier-Stokes equations. To be more specific, we will construct a solution $A$ of \eqref{eq:wholeSpaceLeray}
subjected to the spatial boundary condition 
\begin{align}\label{eq:wholeSpaceLeray2}
	\displaystyle \lim_{|y_0|\to\infty} \int_{B_1(y_0)}|A(y,s)-U_0(y,s)|^2\,dy= 0,
\end{align}
and which is $T$-periodic, i.e.
\begin{align}\label{eq:wholeSpaceLeray3}
  A(\cdot,s)=A(\cdot, s+T),
 \end{align}
for a given $T$-periodic divergence free vector fields $B$ and $U_0$.  Here, $B$ is as above and $U_0$ is defined to be 
\[
U_0(y,s)=\sqrt t e^{t\Delta}a_0 (x).
\]
Thus, $U_0$ serves as the boundary value of the system for $A$ and encodes information about $a_0$ as boundary data.

It can be shown that $U_0$ satisfies the following assumption (see \cite{BT2}):
\begin{assumption} \label{AU_0}
The vector field $U_0(y,s) :\R^3 \times \R \to \R^3$ is continuously differentiable in $y$ and $s$,  periodic in $s$ with period $T>0$, divergence free, and satisfies
\begin{align*}
& 	\partial_s U_0 
-\Delta U_0-\frac 1 2 U_0-\frac 1 2 y\cdot \nabla U_0 = 0, 
\\& U_0\in L^\infty (0,T;L^4\cap L^q(\R^3)), 
\\& \partial_s U_0\in L^\infty(0,T;L_{\mathrm{\loc}}^{6/5}(\R^3)),
\end{align*}
and
\[
\sup_{s\in [0,T]}\|U_0  \|_{L^q(\R^3\setminus B_R)}\leq \Theta(R),
\]
for some $q\in (3,\infty]$ and $\Theta:\R_+\to \R_+$ such that $\Theta(R)\to 0$ as $R\to\infty$.
\end{assumption}

When we say that $A$ solves \eqref{eq:wholeSpaceLeray}--\eqref{eq:wholeSpaceLeray3}, we mean it is a periodic weak solution in the following sense.

\begin{definition}[Periodic weak solution]
\label{def:periodicweaksolutionR3} 
Let $U_0$ satisfy Assumption \ref{AU_0} and assume $B$ is $T$-periodic and divergence free. 
The field $A$ is a periodic weak solution to \eqref{eq:wholeSpaceLeray} in $\R^3\times (0,T)$ if it is divergence free, if 
\begin{equation}\notag
A-U_0\in L^\infty(0,T;L^2(\R^3))\cap L^2(0,T;H^1(\R^3)),
 \end{equation} 
and if
\EQ{\label{u.eq-weak}
&\int_0^T \bigg( (A,\partial_s f)-(\nabla A,\nabla f)
\bigg)  \,ds 
\\& = -\int_0^T \bigg( \frac 1 2 A+\frac 1 2 y\cdot\nabla A-A\cdot\nabla A -A\cdot \nabla B -B\cdot\nabla A,f  \bigg)  \,ds,
} 
holds for all $f \in \mathcal D_T$.  
This latter condition implies that $A(0)=A(T)$.

\end{definition}

In the above, $A-U_0$ is in the energy class. This hints that we will use energy estimates and compactness arguments to construct $A-U_0$.  Attempting to do so, we hit a roadblock when trying to obtain the formal bound
\[
 \int_0^T\int (A\cdot\nabla U_0)\cdot A \,dy\,ds\leq \gamma \|  A\|_{H^1}^2,
\]
where $\gamma$ is a prescribed, small parameter.  Basically, $U_0$ isn't small and we therefore need to replace it with something that is both small and asymptotically comparable to $U_0$. This is possible because $a_0$ has some decay at spatial infinity.  This allowed us to prove the following lemma in \cite{BT1}.

\begin{lemma} 
\label{lemma:W}
Fix $q\in (3,\infty]$ and suppose $U_0$ satisfies Assumption \ref{AU_0} for this $q$. 
Fix $Z\in C^\infty(\R^3)$ with $0 \le Z \le 1$, $Z(x)=1$ for $|x|>2$ and $Z(x)=0$ for $|x|<1$.
For any $\alpha\in (0,1)$, there exists $R_0=R_0(U_0,\alpha)\ge 1$ so that letting $\xi(y) =Z(\frac y{R_0})$ and setting
\begin{equation}
  W (y,s)= \xi(y) U_0(y,s) + w(y,s),
\end{equation}
where 
\begin{equation}
w(y,s)=\int_{\R^3}\nabla_y \frac 1 {4\pi |y-z|} \nabla_z \xi(z) \cdot U_0 (z,s) \,dz,
\end{equation}
we have that $W$ is locally continuously differentiable in $y$ and $s$, $T$-periodic, divergence free,
 $U_0 - W \in L^\infty(0,T; L^2(\R^3))$, and 
\begin{equation}\label{ineq:Wsmall}
\|W\|_{L^\infty(0,T;L^q(\R^3))}\leq \alpha, %
\end{equation} 
\begin{equation}\label{WL4.est}
\norm{W}_{L^\infty(0,T;L^4(\R^3))}\leq c(R_0,U_0),
\end{equation}
and
\begin{equation}
\label{LW.est}
\norm{	\partial_s W
-\Delta W-\frac 1 2 W-\frac 1 2 y\cdot \nabla W}_{L^\infty(0,T; H^{-1}(\R^3))} \leq c(R_0,U_0), %
\end{equation}
where $c(R_0,U_0)$ depends on $R_0$ and quantities associated with $U_0$ which are finite by Assumption \ref{AU_0}.  
\end{lemma}  
For the proof, see \cite{BT1} (a more technical but more robust proof is given in \cite{BT2} which is also valid on the half-space).

It is now clear that we should seek a solution of the form \[U:=A-W.\]
A bit of calculus reveals that the weak formulation for $U$ is: for all $f \in \mathcal V$ and a.e.~$s\in (0,T)$,
\begin{align}\label{perturbed-Leray}
\frac d {ds}(U,f)
&=	
		- (\nabla U,\nabla f) 
		+ (%
		\frac 12 U+\frac 12 y\cdot \nabla U, f) - (U \cdot\nabla U, f)
\\ 
\notag& 
		-(W\cdot\nabla U+U\cdot \nabla W + U\cdot \nabla B +B\cdot \nabla U ,f)
\\ 
\notag& 		-\langle \mathcal{R}(W),f\rangle,
\end{align}
where 
\begin{align}\label{RW.def}
\mathcal{R}(W) & :=  	\partial_s W %
-\Delta W
\\&-\frac 12W-\frac 12y\cdot \nabla W +B\cdot\nabla W+W\cdot\nabla B+ W\cdot\nabla W .\notag
\end{align}
In \cite{BT1}, the solution, which evolves from data in $L^3_w$, satisfies the local energy inequality.  It turns out that the same is true of $a$ in the present context, although this was not pursued in \cite{BT3}.  In order to prove this, we modify the formulation of the problem by mollifying the drift velocity.  To this end, for $\epsilon>0$, let $\eta_\epsilon(y)=\epsilon^{-3}\eta(y/\epsilon)$ where $\eta\in C_0^\infty$ is fixed and satisfies $\int_{\R^3}\eta\,dy=1$. The mollified $U$ equation becomes
\begin{align}\label{eq.mollified-perturbed-Leray}
\frac d {ds}(U,f)
&=	
		- (\nabla U,\nabla f) 
		+ (%
		\frac 12 U+\frac 12 y\cdot \nabla U, f)  
		- ((\eta_\e *U) \cdot\nabla U, f)
\\
\notag &  \quad
		-(W\cdot\nabla U+U\cdot \nabla W + U\cdot \nabla B +B\cdot \nabla U ,f)	
		-\langle \mathcal{R}(W),f\rangle.
\end{align}

\noindent{\bf Step 2: Auxiliary problem -- approximation.}  
We use the Galerkin method as in \cite{BT1}.
Let $\{{a_k}\}_{k\in \N} \subset \mathcal V$ be orthonormal in $L^2_\si (\R^3)$ and its span be dense in $X$.
For a fixed $k$, we look for an approximation solution of the form $U_k(y,s)= \sum_{i=1}^k b_{ki}(s)a_i(y)$.
We first prove the existence of and \emph{a priori} bounds for $T$-periodic solutions $b_k=(b_{k1},\ldots,b_{kk})$ to the system of ODEs
\begin{align}\label{eq:ODE}
\frac d {ds} b_{kj} = & \sum_{i=1}^k A_{ij}b_{ki} +\sum_{i,l=1}^k B_{ilj} b_{ki}b_{kl} +C_j,%
\end{align}
for $j\in \{1,\ldots,k\}$,
where
\begin{align}
\notag A_{ij}&=- (\nabla a_{i},\nabla a_j) 
		+ (%
		\frac 12 a_i+\frac y2 \cdot \nabla a_i, a_j) 
\\\notag &\quad\,	 -(a_i\cdot \nabla (W+B),a_j)
		- ((W+B)\cdot\nabla a_i, a_j)
\\\notag B_{ilj}&=- ((\eta_\e*a_i) \cdot\nabla a_l, a_j)
\\\notag C_j&=-\langle \mathcal R (W),a_j\rangle.
\end{align}

We will show that,
 for any $k\in \mathbb N$ and $\e>0$, the system of ordinary differential equations \eqref{eq:ODE} has a $T$-periodic solution $b_{k}\in H^1(0,T)$.  Letting
\begin{equation} \notag
U_{k,\e}(y,s)=\sum_{i=1}^k b_{ki}(s)a_i(y),
\end{equation}it is possible to show that
\begin{equation}\label{ineq:uniformink}
||U_{k,\e}||_{L^\infty (0,T;L^2(\R^3))} + ||U_{k,\e}||_{L^2(0,T;H^1(\R^3))}<C,
\end{equation}where $C$ is independent of $k$ and $\e$.   We sketch the proof:  First, we solve the problem \eqref{eq:ODE} for an arbitrary $U^0\in \operatorname{span}(a_1,\ldots, a_k)$
in a possibly short time interval $[0,T']$ (this is a classical ODE problem for which smooth solutions exist) and obtain energy estimates for the solution using Gronwall's inequality.  In particular, by applying the estimates from \cite[(2.18)--(2.22)]{BT2} and \cite[(3.11)--(3.16)]{BT3}, we see that
\begin{equation}
\begin{split}
e^{s/4} ||U_{k,\e}(s)||_{L^2}^2
&\leq ||U^{0}||_{L^2}^2 +  \int_0^{  T} e^{\tau/4}  C_2 \,d\tau
\\
& \le  ||U^{0}||_{L^2}^2 + e^{T/4} C_2 T,
\end{split}
\end{equation}
for all $s\in [0, T']$. Using this inequality, we can extend the existence time to $T'=T$ and choose $\rho>0$ (independent of $k$ and $\e$) so that 
\begin{equation}\notag
 ||U^{0}||_{L^2}\leq \rho \Rightarrow ||U_{k,\e}(T')||_{L^2}\leq \rho.
\end{equation}
The mapping $\mathcal T:B_\rho^k\to B_\rho^k$ given by $\mathcal T(b_{k}(0))=b_k(T)$, where $ B_\rho^k$ is the closed ball of radius $\rho$ in $\R^k$, is continuous. Thus $\mathcal T$ has a fixed point by the Brouwer fixed-point theorem, i.e.~there exists some $U^{0}\in \operatorname{span}(a_1,\ldots,a_k)$ so that $b_k(0)=b_k(T)$.  It follows that $U_{k,\e}(T)=U^0(T)=U_{k,\e}(0)$, i.e.~$U_{k,\e}$ is $T$-periodic.
Finally, the uniform bound \eqref{ineq:uniformink} follows from the energy estimate.

\medskip
\noindent{\bf Step 3: Auxiliary problem -- convergence.}  We have two limiting parameters, $\e$ and $k$.  We first take a limit in $k$ to obtain a solution $U_\e$ from the sequence $\{U_{k,\e}\}$. Standard arguments (e.g.~those in \cite{Temam}) imply that, for $T>0$ and for any $\epsilon>0$, there exists $T$-periodic 
$U_\e\in {L^2(0,T;H^1(\R^3))}$
(with norm bounded independently of $\epsilon$) and a subsequence of $\{U_{k,\e}\}$ (still denoted by $U_{k,\e}$) so that 
\begin{align*}
& U_{k,\e}\rightarrow U_\epsilon \mbox{~weakly in}~L^2(0,T;X),
\\& U_{k,\e}\rightarrow U_\epsilon \mbox{~strongly in}~L^2(0,T;L^2(K))  \mbox{~for all compact sets~}K\subset \R^3,
\\& U_{k,\e}(s)\rightarrow U_\epsilon(s) \mbox{~weakly in}~L^2 \mbox{~for all}~s\in [0,T].
\end{align*}
The weak convergence guarantees that $U_\epsilon(0)=U_\epsilon(T)$.  The limit $U_\e$ is a periodic weak solution of the mollified perturbed Leray system.  The variational problem made no mention of the pressure.  Since we want to prove $a$ satisfies the local energy inequality, we need to construct the pressure explicitly.  To do this note that $B\cdot\nabla W +W\cdot\nabla B\in L^\I(H^{-1})$.  Then, the argument of \cite[Proof of Theorem 2.4]{BT1} applies essentially verbatim and we conclude that there exists $p_\e$ so that $U_\e$ and $p_\e$ are a classical solution to the mollified, $W$- and $B$-perturbed, $T$-periodic Leray equations and, furthermore, $p_\e\in L^{5/3}(\R^3\times [0,T])$ with norm bounded uniformly in $\e$.  Because these are classical solutions, they satisfy local energy \emph{equalities}.  

We finally let $\e\to 0$.  
Because $U_\epsilon$ are bounded independently of $\epsilon$ in \[L^\infty (0,T;L^2(\R^3))  \cap   L^2(0,T;H^1(\R^3)),\] and $U_\e$ is a weak solution of \eqref{eq.mollified-perturbed-Leray} with $W$ bounded by Lemma \ref{lemma:W}, there exists a vector field $U\in L^\infty (0,T;L^2(\R^3))\cap L^2(0,T;H^1(\R^3))$ and a sequence $\{ U_{\epsilon_k}\}$ of elements of $\{U_\epsilon \}$ so that
\begin{align*}
& U_{\epsilon_k} \rightarrow U \mbox{~weakly in}~L^2(0,T;X)
\\& U_{\epsilon_k}\rightarrow U \mbox{~strongly in}~L^2(0,T;H(K)) ~ \forall \mbox{~compact sets $K\subset \R^3$}
\\& U_{\epsilon_k}(s)\rightarrow U(s) \mbox{~weakly in}~L^2 \mbox{~for all}~s\in [0,T],
\end{align*}
as $\epsilon_k\to 0$.   Furthermore, since $p_{\epsilon_k}$ are uniformly bounded in $L^{5/3}(\R^3\times [0,T])$ we can extract a subsequence (still denoted $p_{\epsilon_k}$) so that 
\begin{equation}
p_{\epsilon_k}\rightarrow p_U \mbox{~weakly in}~L^{5/3}(\R^3\times [0,T]),
\end{equation}
for some distribution $p_U\in L^{5/3}(\R^3\times [0,T])$ and this convergence is strong enough to ensure that $(U,p_U)$ solves \eqref{eq:wholeSpaceLeray} in the distributional sense.

We now prove that $A=U+W$ satisfies a local energy inequality.  Earlier we noted that $A_\e=U_\e+W$ satisfies a local energy equality. Indeed, if $\psi \in \mathcal C_0^\I(\R^4)$ is non-negative, then,
multiplying the $A_\e$-equation by $2A_\e \psi$ and
integrating over $\R^3 \times \R$,
\begin{align*}
\iint &\bke{ 2|\nabla A_\e|^2+\tfrac 12 |A_\e|^2} \psi\,dy\,ds 
\\
&= \iint |A_\e|^2(\partial_s \psi +\Delta \psi)\,dy\,ds 
\\
&+\iint \bke{ |A_\e|^2 \bkt{ \eta_\e * U_\e+W +B-\tfrac 12y } + 2p_{A_\e} A_\e} \cdot \nabla \psi\,dy\,ds
\\&+ \iint \bkt{ (\eta_\e * U_\e - U_\e)\cdot \nb W - A_\e \cdot\nabla B}\cdot 2A_\e\psi \,dy\,ds,
\end{align*}
where 
\[
p_{A_\e} = \sum_{i,j} R_iR_j \big \{
(\eta_\e * U_{\e,i}+W_i)U_{\e,j}
+A_{\e,i}(W_{j}
+ B_j)+B_i A_{\e,j}\big \}.
\]

As in the discussion following \cite[(A.51)]{CKN}, the left hand side of the above equality is lower semi-continuous as $\e\to 0$, while each term on the right-hand side converges to the corresponding $A$ term.   Therefore, after several cancellations, $A$ satisfies  
\begin{align}\label{ineq.localA}
&\iint \bke{ 2|\nabla A|^2+\tfrac 12 |A|^2} \psi\,dy\,ds 
\\\notag &\leq \iint |A|^2(\partial_s \psi +\Delta \psi)\,dy\,ds 
\\\notag
&+\iint \bke{ |A|^2 \bkt{ A +B-\tfrac 12y} + 2p_A A} \cdot \nabla \psi\,dy\,ds
\\\notag &- \iint \bkt{ A \cdot\nabla B}\cdot 2A\psi \,dy\,ds,
\end{align}
where 
\[
p_A = \sum_{i,j} R_iR_j \big \{
A_iA_{j} 
+A_{i}B_j+B_i A_j\big \}.
\]

\medskip

\noindent{\bf Step 4: Loose ends}

At this point we have constructed a solution $A$ and associated pressure $p_A$ to \eqref{eq:wholeSpaceLeray}.  Since $B$ is prescribed, this gives us a solution $u=A+B$ of \eqref{eq:Leray-rot}, and this corresponds to a solution $v=a+b$ of the Navier-Stokes equations.  Except for the local energy inequality \eqref{ineq.localv}, all the conclusions of the theorem follow identically to \cite[Proof of Theorem 1.4]{BT3}. 

For the local energy inequality \eqref{ineq.localv}, we first show that \eqref{ineq.localA} for $A,p_A$ is equivalent to 
the following $b$-perturbed local energy inequality for $a$ and $p_a$: if $\phi\in C^\I_0(\R^3\times \R^3_+)$ is non-negative then
\EQ{ \label{ineq.locala}
&2\iint |\nabla a|^2\phi\,dx\,dt 
\leq  \iint |a|^2(\partial_t \phi + \Delta\phi )\,dx\,dt %
\\&\quad + \iint [|a|^2 (a+b) + 2p_a a]
\cdot\nabla \phi\,dx\,dt - 2\iint (a\cdot\nabla b) \cdot (a\phi)\,dx\,dt.
}
To see this, let $\phi(x,t)=\frac 1 {\sqrt t} \psi(y,s)$, where $\psi$ is a test function from \eqref{ineq.localA}. Using the relationships between $a$, $A$, $b$, $B$, $p_a$, $p_A$, $y$, $x$, $s$, and $t$, it is possible (after a lengthy computation) to derive \eqref{ineq.locala} directly from \eqref{ineq.localA}.  The other direction follows from the same computation.

To show \eqref{ineq.localv}, note that $b$ is smooth for $t>0$, and the identities from integration by parts of
\[
\iint (\pd_t a - \De a + (a+b)\cdot\nabla  a + a \cdot \nabla b + \nabla p_a)\cdot 2b \phi=0
\]
and
\[
\iint (\pd_t b - \De b + b\cdot \nabla b+\nabla p_b)\cdot 2v \phi=0
\]
are valid. The sum of these identities and \eqref{ineq.locala} gives  \eqref{ineq.localv}. 

This completes the proof of Theorem \ref{thrm.new} in the DSS case.

\subsection{Proof of Theorem \ref{thrm.new}, SS case}

The SS case can be proved directly by constructing stationary  solutions to \eqref{eq:stationaryLeray}.  Alternatively, they can also be obtained as a limit of DSS solutions.  We give a sketch of the proof of Theorem \ref{thrm.new} using the limit approach since it illustrates the usefulness of the a priori bound \eqref{ineq.apriorilocal} satisfied by local Leray solutions.  This approach was also used in \cite[Proof of Theorem 1.3]{BT1} to obtain SS solutions as a limit of DSS solutions.

\begin{proof}%

Assume $v_0$ is SS and belongs to $\Bp$.  Let $a_0$ and $b_0$ come from Lemma \ref{lemma.profileslicing} for a small constant $\e_0$ (smaller than the Koch-Tataru constant) and let $b$ be the unique mild solution to the Navier-Stokes equations \eqref{eq:NSE} with initial data $b_0$.  Let $a_0^{(k)}=a_0$ where $k\in \N$.  Clearly, $a_0^{(k)}$ is $\la_k$-DSS with $\la_k=2^{2^{-k}}$.  Thus, by the DSS case, there exists a $\la_k$-DSS solution $v^{(k)}=a^{(k)}+b$ to the Navier-Stokes equations.  Moreover, each $a^{(k)}$  satisfies \eqref{ineq.apriorilocal} but where the constant $\si$ is modified (in comparison to \eqref{ineq.apriorilocalconstant}) to also depend on $b$ (we omit a proof of this but it follows from the fact that $a$ satisfies a linear perturbation of the Navier-Stokes equations and the argument from \cite{JiaSverak}).  Since $b$ is independent of $k$, this gives a uniform bound on the sequence $a^{(k)}$ and allows us to apply a convergence argument to obtain a SS solution to the Navier-Stokes equations satisfying the properties listed in the statement of Theorem \ref{thrm.new} -- see \cite[Proof of Theorem 1.3]{BT1} for details.
\end{proof}

\section*{Acknowledgments}
The research of both authors was partially supported by the NSERC grant 261356-13 (Canada). That of Z.B. was also partially supported by the NSERC grant 251124-12.

\bigskip

Zachary Bradshaw, Department of Mathematics, University of Arkansas, Fayetteville, AR 72701, USA;
e-mail: 
zb002@uark.edu

\medskip

Tai-Peng Tsai, Department of Mathematics, University of British
Columbia, Vancouver, BC V6T 1Z2, Canada;
e-mail: ttsai@math.ubc.ca

\end{document}